\documentclass[12pt]{amsart}
\usepackage{amsfonts} 
\usepackage[all]{xy}
\usepackage{graphicx}
\usepackage{amssymb}
\usepackage{color}
\makeatletter
\newdimen\sectionruledimen
\def\makeparrule{
  \def\par{
   \endgraf\nobreak\vskip\lineskip\nointerlineskip
   \hbox to\hsize{\hskip\sectionbackskip\leaders\hrule height
   \sectionruledimen\hfil}
   }
 }

\def\section{%
 \@startsection {section}{1}{\sectionbackskip}{-10pt plus -1.2ex minus -.2ex}%
  {.5pt}{\normalsize\bf\makeparrule}%
 }%
\newdimen\sectionbackskip
\sectionbackskip\z@
\newdimen\sectionruledimen
\sectionruledimen\z@

\def\punteada{\leaders\hbox{$\m@th \mkern1.5mu - \mkern1.5mu$}\hfill}
\makeatother

\newtheorem{theorem}{Theorem}
\newenvironment{theor}{\smallskip\begin{trivlist}
   \item[\hspace{\labelsep}{\noindent\bf Theorem.}]\it
   }{\end{trivlist}\smallskip}

\newenvironment{lemma}{\smallskip\begin{trivlist}
   \item[\hspace{\labelsep}{\noindent\bf Lemma.}]\it
   }{\end{trivlist}\smallskip}

\newenvironment{propo}{\smallskip\begin{trivlist}
   \item[\hspace{\labelsep}{\noindent\bf Proposition.}]\it
   }{\end{trivlist}\smallskip}

\newcommand{\cuadro}{\hfill{$\qed$}}

\def\lftcol#1{\vbox {\halign {##\hfil \cr #1\cr }}}
\long\def\direc#1#2{\hbox to \hsize{
     \lftcol{\hsize=7cm#1}\hfill\lftcol{\hsize=7cm#2}}}
\newcommand{\A}{\mathbb{A}}

\newcommand{\D}{\mathbb{D}}
\newcommand{\E}{\mathbb{E}}

\newcommand{\N}{\mathbb{N}}

\newcommand{\R}{\mathbb{R}}
\newcommand{\s}{\mathbb{S}}
\newcommand{\T}{\mathbb{T}}
\newcommand{\x}{\mathbb{X}}
\newcommand{\X}{\mathbb{X}}
\newcommand{\Z}{\mathbb{Z}}

\newcommand{\mH}{\mathcal{H}}

\newcommand{\coh}{\hbox{\rm coh}}
\renewcommand{\deg}{\hbox{\rm deg}}

\newcommand{\Der}{\hbox{\rm Der}}

\renewcommand{\dim}{\hbox{\rm dim}}

\newcommand{\Ext}{\hbox{\rm Ext}}
\newcommand{\End}{\hbox{\rm End}}

\newcommand{\Hom}{\hbox{\rm Hom}}

\renewcommand{\mod}{\hbox{\rm mod}}

\newcommand{\Root}{\hbox{\rm Root}}
\newcommand{\Spec}{\hbox{\rm Spec}}

\newcommand{\bulito}{\ {\scriptstyle \bullet}\, \ }

\def\subsetnoteq{\mathbin{\hbox{$\subseteq \joinrel \hskip-8pt \lower3pt
                 \hbox{$\scriptscriptstyle /$}\ $}}}

\def\raya{\raise1.5pt\hbox to 25pt{\vrule height1.5pt depth-1pt
           width25pt}}
\def\rayita{\raise2pt\hbox to 7.5pt{\vrule height1.5pt depth-1pt
           width7.5pt}}

\begin{document}

\title[Mahler measure of Coxeter polynomials ]{On the Mahler measure of the Coxeter polynomial of  algebras.}
\author[Jos\'e A. de la Pe\~na]{Jos\'e A. de la Pe\~na}

\begin{abstract}  
Let $A$ be a finite dimensional algebra over an algebraically closed field $k$. Assume $A$ is a basic connected and triangular algebra with $n$ pairwise non-isomorphic simple modules. We consider the {\em Coxeter transformation} $\phi_A(T)$ as the automorphism of the Grothendieck group $K_0(A)$ induced by the Auslander-Reiten translation $\tau$ in the derived category $\Der^b(\mod_A)$ of the module category $\mod_A$ of finite dimensional left $A$-modules. We say that $A$ is of {\em cyclotomic type} if the characteristic polynomial $\chi_A$ of $\phi_A$ is a product of cyclotomic polynomials, equivalently, if the {\em Mahler measure} $M(\chi_A)=1$. In \cite{Pe} we have considered the many examples
of algebras of cyclotomic type in the representation theory literature. In this paper we study the Mahler measure of the Coxeter polynomial of {\em accessible algebras}. In 1933, D. H. Lehmer  found that the polynomial $T^{10} + T^9 - T^7 - T^6 - T^5 - T^4 - T^3 + T + 1$ has Mahler measure $\mu_0 = 1.176280 . . .$, and he asked if there exist any smaller values  exceeding 1. In this paper we prove that for any accessible algebra $A$ either $M(\chi_A)=1$ or $M(\chi_B) \ge \mu_0$ for some convex subcategory $B$ of $A$. We introduce {\em interlaced tower of algebras} $A_m,\ldots,A_n$ with $m \le n-2$ satisfying
$$\chi_{A_{s+1}} =(T+1) \chi_{A_s} - T \chi_{A_{s-1}}$$
for $m+1 \le s \le n-1$. We prove that, if ${\rm Spec \,}\phi_{A_n} \subset \s^1 \cup \R^+$ and $A_n$ is not 
of cyclotomic type then  $M(\chi_{A_m}) <M(\chi_{A_n})$. 
\noindent
\end{abstract}
\maketitle

\section{Introduction.}
Assume throughout the paper that $k$ is an algebraically closed field.
We assume that $A$ is a triangular finite dimensional $k$-algebra, in particular, $A$ has finite global dimension. 
In this situation, the bounded derived category $\Der(A)=\Der^b(\mod_A)$ has Serre duality of the form
$$
\Hom(X,Y[1])=\Hom(Y,\tau X)
$$
where $\tau$ is a self-equivalence which serves as the Auslander-Reiten translation of $\Der(A)$. In this setting the {\em Grothendieck group} $K_0(A)$ of the category $\mod_A$ of finite dimensional (right) $A$-modules, formed with respect to short exact sequences, is naturally isomorphic to the Grothendieck group of the derived category, formed with respect to exact triangles. 
The Coxeter transformation $\phi_A$ is the automorphism of the Grothendieck group $K_0(A)$ induced by the Auslander-Reiten translation $\tau$ in the derived category $\Der(A)$. The characteristic polynomial $\chi_A(T)$ of $\phi_A$ is called the {\em Coxeter polynomial} $\chi_A(T)$ of $A$, or simply $\chi_A$. It is a monic self-reciprocal polynomial, therefore it is $\chi_A(T)= a_0 + a_1 T + a_2 T^2 + \ldots + a_{n-2} T^{n-2} + a_{n-1}T^{n-1} + a_n T^n \in \Z[T]$, with
$a_i=a_n-i$ for $0 \le i \le n$, and $a_0=1=a_n$. Sometimes, for the sake of clarity, we write $a_i(A)$ for these coefficients.  

Consider the roots $\lambda_1(A),\ldots,\lambda_n(A)$ of $\chi_A$, the so called  {\em spectrum} of $A$. There is a number of measures associated to the absolute values
$|\lambda|$ for $\lambda$ in the spectrum ${\rm Spec}(\phi_A)$ of $A$. Three important measures are the following:
the {\em spectral radius} of $A$ is defined as $\rho_A = {\rm max \, } \{|\lambda| : \lambda \in {\rm Spec}(\phi_A) \}$,
the {\em Mahler measure} of $\chi_A$ is $M(\chi_A)= \prod_{i=1}^n  |\lambda_i|$ and the {\em energy function} 
of $\phi_A$ is $e(\phi_A)=\sum\limits_{i=1}^n |\lambda_i(A)|$. In a recent paper \cite{Pe} we show that  always $e(\phi_A) \ge n$ and equality is reached when 
$A$ is of {\em cyclotomic type}, that is, $\chi_A(T)$ factorizes as product of cyclotomic polynomials (equivalently  $\rho_A =1$ or $M(\phi_A) = 1$).
Many finite dimensional algebras are known to be of cyclotomic type: hereditary algebras of finite or tame representation type, canonical algebras, some extended canonical algebras and many others, see \cite{LePe3,Pe} and section 2 for more examples.

By a celebrated result of Kronecker \cite{Kr}, a monic integral polynomial $p$ has $M(p)=1$ if and only if $p$ factorizes as product of cyclotomic polynomials. 
Rather little is known, however, about values of the Mahler measure near 1. In 1933, D. H. Lehmer  found that
the polynomial
$$T^{10} + T^9 - T^7 - T^6 - T^5 - T^4 - T^3 + T + 1$$
\noindent
has Mahler measure $\mu_0 = 1.176280 . . .$, and he asked if there exist any smaller values  exceeding 1. In fact, the polynomial above is the Coxeter polynomial
of the hereditary algebra whose underlying graph $[2,3,7]$ is depicted below:
$$
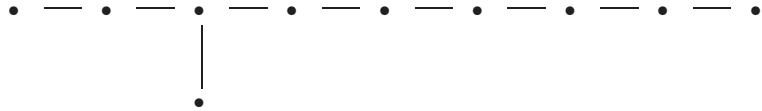
\xymatrix@C14pt{\bulito \ar@{-}[r]&\bulito \ar@{-}[r]&\bulito \ar@{-}[r] \ar@{-}[d]&\bulito \ar@{-}[r]&\bulito \ar@{-}[r]&\bulito \ar@{-}[r]&\bulito \ar@{-}[r]&\bulito \ar@{-}[r]&\bulito \\&&\bulito&&&&&& }$$
 \noindent 
Mossingshoff's web site provides many examples of polynomials $p$ with small Mahler measure $M(p)>1$. See \cite{Mo1,Mo2}.
The purpose of this work is to consider the largely unexplored topic of the Mahler measure $M(\chi_A)$ of the Coxeter polynomial of an algebra $A$, 
the many examples of algebras with Mahler measure 1 and other examples serving as motivation. In Section 4 we shall prove the following theorem, which is to be understood as
a representation theory variant of Lehmer's problem.
\vskip.4cm
\noindent
{\bf Theorem 1} {\em Let $A$ be an accessible algebra not of cyclotomic type. Then there is a convex subcategory $B$ of $A$ satisfying the following properties:

(a) $B$ is minimal not of cyclotomic type, that is, if $C$ is any proper convex subcategory of $B$,  then $C$ is of cyclotomic type;

(b) the Mahler measure of $B$ is $M(\chi_B) \ge \mu_0$.}
\noindent
\vskip.5cm

Section 2 recalls some general features concerning algebras of cyclotomic type. Section 3 presents some simple results on the Mahler measure. The most precise considerations and examples will be worked for the class of accessible algebras. Introduced in \cite{LePe5}, we say that an algebra $A$ is {\em accessible from} $B$ if there is a sequence $B=B_1,B_2,\ldots,B_s=A$ of algebras such that each $B_{i+1}$ is a one-point extension (resp.\ coextension) of $B_i$ for some exceptional $B_i$-module $M_i$.
As a special case, a $k$-algebra $A$ is called {\em accessible} if $A$ is accessible from $k$. In order to prove Theorem 1, we shall provide a complete classification, up to derived equivalence, of the minimal algebras not of cyclotomic type. 

In Section 5 we consider towers of algebras. {\em Interlaced tower of algebras} are certain sequences of triangular algebras $A_m,\ldots,A_n$, 
such that $A_s$ is a basic algebra with $s$ simple modules and, among others, satisfy the common feature that
$$\chi_{A_{s+1}} =(T+1) \chi_{A_s} - T \chi_{A_{s-1}}$$
for $s=m+1,\ldots,n-1$. Examples include sequences of nested hereditary subcategories, sequences of canonical algebras $C_{(p_1(s) ,p_2(s) ,p_3(s))}$ with
$p(s) \le p(s+1)$, sequences of extended canonical algebras $E_{\langle p_1(s), p_2(s), p_3(s) \rangle}$ with $p(s) \le p(s+1)$ and some {\em accessible towers of algebras} where $A_{s+1}$ is a one-point extension (or coextension) of an accessible algebra $A_s$ by an exceptional $A_s$- module $M_s$ such that the perpendicular category $M_s^\perp$ formed in the derived category $\Der(A_s)$ is triangular equivalent to $\mod (A_{s-1})$, for $s=m+1,\ldots,n-1$, 
see \cite{LePe2,LePe3,LePe4,LePe5,Pe} for related concepts.
 \vskip.4cm
{\bf Theorem 2}
{\em Consider an interlaced tower of algebras $A_m,\ldots,A_n$ with $m \le n-2$. If ${\rm Spec \,}\phi_{A_n} \subset \s^1 \cup \R^+$ then  either all $A_i$ are of cyclotomic type or
$M(\chi_{A_m}) <M(\chi_{A_n})$. In the later case, $M(\chi_{A_n})<M(\chi_{A_m})\prod_{s=m+1}^{n-1} L(\chi_{A_s})$, where $L(p)$ denotes the length of a 
polynomial $p$}.
\noindent
\vskip.5cm
Our results make use of well established techniques in the {\em representation theory of algebras} but only reduced use of results from the {\em theory of polynomials} and {\em trascendental number theory}, where Mahler measure has its natural habitat. In forthcoming work we shall extend our investigations on the largely unexplored overlapping area of these important subjects.  

The results of this paper were obtained during a sabbatical stay at CIMAT, Guanajuato. They were presented at seminars in M\'exico, Morelia, Guanajuato, Bielefeld and Tor\'un 
during 2013.

\section{The spectrum of algebras of cyclotomic type} 

\subsection{}
If the spectrum of $A$ lies in the unit disk, then all roots of $\chi_A(T)$ lie on the unit circle, hence $A$ has spectral radius $\rho_A=1$. Clearly, for fixed degree there are only finitely many monic integral polynomials with this property.  The following result of Kronecker is essential in the discussion, \cite{Kr}.

\begin{theor} {\rm [Kronecker]}
Let $f$ be a monic integral polynomial whose spectrum is contained in the unit disk. Then all roots of $f$ are roots of unity or $0$. Equivalently, $f$ factors into cyclotomic polynomials and $T^m$, for some integer $m \ge 0$.
\end{theor}

We recall some facts on {\em cyclotomic polynomials}.

The $n$-cyclotomic polynomial $\Phi_n(T)$ is inductively defined by
the formula
$${\textstyle T^n-1=\prod\limits_{d\mid n}\Phi_d(\mathcal{T}).}$$

The {\em M\"obius function\/} is defined as follows:
$$\mu (n)=\begin{cases}
0 &\hbox{if $n$ is divisible by a square}\\
(-1)^r &\hbox{if $n=p_1,\ldots p_r$ is a factorization into distinct
primes.}\end{cases}$$ 
\noindent 
A more explicit expression for the cyclotomic polynomials is given by:

\begin{lemma}
For each $n\ge 2$, we have
$${\textstyle \Phi_n(T)=\prod\limits_{1\le d<n\atop d\mid
n}v_{n/d}(T)^{\mu (d)}}$$
  \vglue-34pt\cuadro
\end{lemma}

In the lemma, we set  $v_n=1+T+T^2+\ldots+T^{n-1}$. Note that $v_n$ has degree $n-1$.

\subsection{}
 A {\em path algebra} $k\Delta$ is said to be of {\em Dynkin type}
if the underlying graph $|\Delta|$ of $\Delta$ is one of the {\em ADE-series}, that is, of type, $\A_n, \D_n$, for some $n \ge 1$ or $\E_p$, for $p=6,7,8$.
The corresponding factorization of the Coxeter polynomial $\chi_{k\Delta}$ is as follows.
\small
\begin{center}
\renewcommand\arraystretch{1.4}
\begin{tabular}{|l|c|c|c|c|} \hline
  Dynkin  type& star symbol & cyclotomic factorization& Coxeter number\\ \hline
  $\A_{n}$& $[n]$  & $\displaystyle\prod_{d|n,d>1}\Phi_d$ & $n+1$\\ \hline
  $\D_{n}$&$[2,2,n-2]$ & $\displaystyle\Phi_2\,\prod_{\frac{d|2(n-1)}{d\neq 1,d\neq n-1}}\Phi_{d}$ & $2(n-1)$ \\\hline
  $\E_{6}$ &$[2,3,3]$ & $\Phi_3\Phi_{12}$& 12 \\ \hline
  $\E_{7}$&$[2,3,4]$& $\Phi_2\Phi_{18}$& 18 \\\hline
  $\E_{8}$& $[2,3,5]$ &$\Phi_{30}$&30\\\hline
\end{tabular}
\end{center}
\normalsize

\subsection{}
There are various instances where an explicit formula for the Coxeter polynomial is known.

Let $A$ be the path algebra of a hereditary star $[p_1,\ldots,p_t]$
with respect to the standard orientation, see
 $$
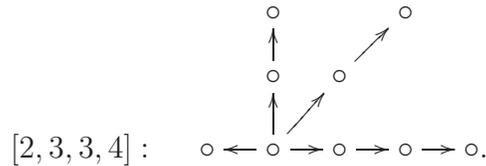
\def\c{\circ}
\xymatrix@C12pt@R12pt{
        &&\c                          &        &\c  &\\
        &&\c\ar[u]                    &\c\ar[ru]&     &\\
[2,3,3,4]:&\c      &\c\ar[l]\ar[u]\ar[r]\ar[ru]&\c\ar[r]&\c\ar[r]&\c.\\
}
$$
Since the Coxeter polynomial $\chi_{A}$ does not depend on the orientation of $A$ we will denote it by $\chi_{[p_1,\ldots,p_t]}$. It follows that
\begin{equation}\label{eq:star_formula}
\chi_{[p_1,\ldots,p_t]}=\prod_{i=1}^{t}v_{p_i}\left((T+1)-T\,\sum_{i^=1}^t \frac{v_{p_i-1}}{v_{p_i}} \right).
\end{equation}
In particular, we have an explicit formula for the sum of coefficients of $\chi_{[p_1,\ldots,p_t]}$ as follows:
$$
\sum\limits_{i=0}^n a_i=  \chi_{[p_1,\ldots,p_t]}(1)=\prod_{i=1}^t p_i\left( 2-\sum_{i=1}^t (1-\frac{1}{p_i})\right).
$$

The next table displays the cyclotomic factorization of extended Dynkin diagrams.
\bigskip
\begin{center}
\renewcommand\arraystretch{1.3}
\begin{tabular}{|c|c|c|c|}\hline
extended Dynkin type & star symbol & weight symbol & Coxeter polynomial \\ \hline
${\tilde \A}_{p,q}$           &     ---     & $(p,q)$         & $(T-1)^2v_p\,v_q$ \\ \hline
${\tilde \D}_{n}$, $n\geq4$   & [2,2,n-2]   & $(2,2,n-2)$   & $(T-1)^2v_2^2v_{n-2}$\\\hline
${\tilde \E}_{6}$& $[3,3,3]$ & $(2,3,3)$    & $(T-1)^2 v_2 v_3^2$ \\\hline
${\tilde \E}_{7}$& $[2,4,4]$ & $(2,3,4)$    & $(T-1)^2 v_2v_3v_4$ \\\hline
${\tilde \E}_{8}$& $[2,3,6]$ & $(2,3,5)$    & $(T-1)^2 v_2v_3v_5$\\\hline
\end{tabular}
\end{center}

Observe taht the Coxeter polynomial of a connected, representation-finite or tame hereditary $k$-algebra $A$ ($k$ algebraically closed) determines the algebra $A$ up to derived equivalence.

\subsection{}
For an algebra $A$ and a left $A$-module $M$ we call 
$$A[M]=\left[\begin{array}{cc}
A & 0\\ M &k
\end{array}\right]$$
\noindent
the {\em one-point extension} of $A$ by $M$. This construction provides an order to deal with {\em triangular algebras}, that is, algebras $kQ/I$, where $I$ is an
ideal of the path algebra $kQ$ for $Q$ a quiver without oriented cycles.

For $A=B[M]$, Happel's long exact sequence \cite{Ha0} relates the Hochschild cohomology groups $H^i(A)$ and $H^i(B)$ in the following way:
$$0\!\to\! H^0(A)\!\to\! H^0(B)\!\to\! \End_B(M)/k\!\to\! H^1(A)\!\to\! H^1(B)\!\to\! \Ext^1_B(M,M)\!\to\! H^2(A)\!\to\! \cdots$$

Therefore, if $M$ is exceptional we get $H^i(A)=H^i(B)$, for $i\ge 0$, and moreover, the cohomology rings $H^*(A)$ and $H^*(B)$ are isomorphic. More general, if $A$ is accessible from $B$,
that is, there exists a chain of algebras $A_1=B, A_2,\ldots, A_s=A$ such that $A_{i+1}$ is a one-point extension or coextension of $A_i$ by an exceptional $A_i$-module $M_i$, then $H^*(A)\cong H^*(B)$ as rings. In particular, if $A$ is accessible (from $k$) then $H^i(A)=0$ for $i>0$ and $H^0(A) =k$.

\subsection{}
Let $c$ be the real root of the polynomial $T^3 - T - 1$, approximately $c=1.325$. As observed in \cite{Xi}, a wild hereditary algebra $A$ associated to a graph $\Delta$ without multiple arrows has $\rho_A > c$ unless $\Delta$ is one of the following graphs:
$${\footnotesize \xymatrix@C12pt@R12pt{ [2,4,5]:&&&&&&& \\1 \ar@{-}[r] &2 \ar@{-}[r]&3 \ar@{-}[r]&4 \ar@{-}[r] \ar@{-}[d]&5 \ar@{-}[r]&6 \ar@{-}[r]&7 \ar@{-}[r]&8  \\&&&9&&&& } \hskip.4cm
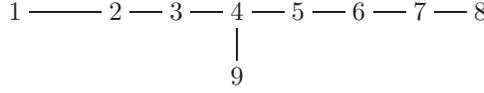
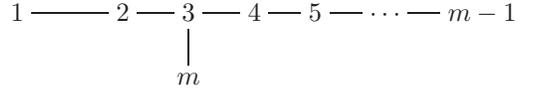
\xymatrix@C12pt@R12pt{[2,3,m]: &&&&&&\\1 \ar@{-}[r]&2 \ar@{-}[r]&3 \ar@{-}[r] \ar@{-}[d]&4 \ar@{-}[r]&5 \ar@{-}[r]&\cdots \ar@{-}[r]&m-1 \\&&m&&&& }}$$
\noindent
In these cases, for $m\ge 8$,
$$c > \rho_{[2,4,5]} > \rho_{[2,3,m]} > \rho_{[2,3,7]}= \mu_0$$
\noindent
where $\mu_0 = 1.176280 . . .$ is the real root of the Coxeter polynomial
$$T^{10} + T^9 - T^7 - T^6 - T^5 - T^4 - T^3 + T + 1$$
\noindent
associated to any hereditary algebra whose underlying graph is $[2,3,7]$. Observe that in these cases, the Mahler measure of the algebra equals the spectral radius.

\subsection{Canonical algebras}

A canonical algebra $\Lambda$ is determined by a weight sequence ${\bf p}=(p_1,\ldots,p_t)$ of $t$ integers $p_i\geq2$ and a parameter sequence ${\bf \lambda}=(\lambda_3,\ldots,\lambda_t)$ consisting of $t-2$ pairwise distinct non-zero scalars  from the base field $k$. (We may assume $\lambda=1$ such that for $t\le 3 $ no parameters occur).  Then the algebra $\Lambda$ is defined by the quiver
$$
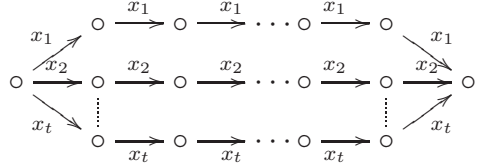
\def\c{\circ}
\xymatrix@C18pt@R10pt{
                                        &\c\ar[r]^{x_1}          & \c\ar[r]^{x_1}         &\cdots\c\ar[r]^{x_1}&\c\ar[dr]^{x_1}&\\
\c\ar[r]^{x_2}\ar[ru]^{x_1}\ar[rd]_{x_t}&\c\ar[r]^{x_2}\ar@{.}[d]&\c\ar[r]^{x_2}          &\cdots\c\ar[r]^{x_2}&\c\ar[r]^{x_2}\ar@{.}[d] &\c    \\
                                        &\c\ar[r]_{x_t}          &\c\ar[r]_{x_t}          &\cdots\c\ar[r]_{x_t}&\c\ar[ru]_{x_t}&    \\
}
$$
satisfying the $t-2$ equations:
$$x_i^{p_i}=x_1^{p_1}-\lambda_i x_2^{p_2},\qquad i=3,\ldots,t.$$
For more than two weights, canonical algebras are not hereditary. Instead, their representation theory is determined by a hereditary category, the category $\coh\X$ of coherent sheaves on a {\em weighted projective line} $\X$, naturally attached to $\Lambda$, see~\cite{GeLe}.
\begin{propo}
  Let $\Lambda$ be a canonical algebra. Then $\Lambda$ is the endomorphism ring of a tilting object in the category $\coh\X$ of coherent sheaves on the weighted projective line $\X$. The category $\coh\X$ is hereditary and satisfies Serre duality in the form $D \Ext^1(X,Y)=\Hom(Y,{\tau X})$ for a self-equivalence $\tau$ which serves as the Auslander-Reiten translation.~\hfill\qed
\end{propo}

Let $\Lambda$ be a canonical algebra with weight and parameter data $({\bf p},{\bf \lambda})$. Then the Coxeter polynomial of $\Lambda$ is given by
$$
\chi_{\Lambda}=(x-1)^2\prod_{i=1}^t v_{p_i}.
$$
\noindent
In particular, canonical algebras are cyclotomic.

A finite dimensional algebra isomorphic to the endomorphism algebra of a tilting object in a (connected) hereditary abelian Hom-finite $k$-category $\mH$ is called a {\em tilted algebra}. By a result of Happel~\cite{Ha1} each such category is derived equivalent to the module category $\mod{H}$ over a hereditary algebra or to the category $\coh\X$ of coherent sheaves on a weighted projective line. The Coxeter polynomials of tilted algebras are therefore the Coxeter polynomials of hereditary or canonical algebras.

\subsection{One-point extensions of canonical algebras} 

Let $\Lambda$ be a canonical algebra of weight type $(p_1,\ldots,p_t)$. 
Then $\Lambda$ is derived equivalent to the category of coherent sheaves $\coh \X$ over a weighted projective line $\X$
of weight type $(p_1,\ldots,p_t)$. An indecomposable $\Lambda$-module $M$ can be chosen as either a vector bundle or a finite length module. 
In case $M$ is a vector bundle of rank one, we call the one-point extension $\Lambda[M]$ an {\em extended canonical algebra}.
If $M$ is a finite dimensional regular module we distinguish two situations:

$(a)$ If $M$ is regular simple in the $i$-th exceptional tube $\T_i$ (of $\tau$-period $p_i$), then the one-point extension $\Lambda[M]$ is tilting-equivalent to the canonical algebra of weight type $(p_1,\ldots,p_i+1,\ldots,p_t)$ having the same parameter sequence as $\Lambda$.

$(b)$ If $M$ has quasi-length $s$ in $\T_i$ (recall this means that $s<p_i$), then $\Lambda[M]$ is derived equivalent to a supercanonical algebra in the sense of \cite{LePe4}, where the linear arms of the canonical algebra with index different from $i$ are kept, and the $i$-th linear arm is changed to the poset $K(p_i,s)$
\footnotesize$$
\xymatrix@R8pt@C14pt{
       &       &             &                   & \star         &         &    \\
1\ar[r]&2\ar[r]&\cdots\ar[r] &p_i-s\ar[r]\ar[ru] & p_i-s+1\ar[r] &\cdots\ar[r]&p_i-1
}
$$\normalsize
We write $\Lambda(i,s)$ for $\Lambda[M]$ and call it a {\em supercanonical algebra of restricted type}.

The structure of the bounded derived category of an extended canonical algebra
$\langle{p_1,\dots,p_t}\rangle$ sensibly depends on the sign of the Euler
characteristic $\chi_{\x}=2-\sum_{i=1}^t(1-1/p_i)$ of the weighted projective line $\x$
associated to $\Lambda$. According to \cite{LePe3}, the description of the derived category of an extended canonical algebra yields a trichotomy.

$(i)$ Positive Euler characteristic: Let $\Lambda$ be a canonical algebra of domestic representation type $(p_1,p_2,p_3)$, and $\Delta$ be the Dynkin diagram $[p_1,p_2,p_3]$. Then the extended canonical algebra $\langle{p_1,p_2,p_3} \rangle$ is derived equivalent to the (wild) path algebra  of a quiver $Q$ having extended Dynkin type $\tilde{\tilde{\Delta}}$.

$(ii)$ Euler characteristic zero: Consider a canonical algebra with a tubular weight sequence $(p_1,\ldots,p_t)$, we shall assume that $2\le p_1\le p_2\le \cdots \le p_t$. Then the extended canonical algebra
$\langle {p_1,\dots,p_t} \rangle$ is derived canonical of type $(p_1,\ldots,p_{t-1},p_t+1)$.

$(iii)$ Negative Euler characteristic: Let $\x$ be a weighted projective line of negative Euler characteristic, let $A= \langle {p_1,\dots,p_t} \rangle$ be the corresponding extended canonical algebra and $R$ be the $\Z$-graded singularity attached to $\x$. Then the derived category $\Der(\mod A)$ is triangle-equivalent to the triangulated category of graded singularities $D_{{\rm sing}}^{{\Z}}(R)$ of $R$, where the superscript $\Z$ refers to the grading.

\subsection{} In \cite{Pe}, we showed that several of the common properties of classes of algebras presented in the above sections are consequence of the cyclotomic condition. 
In particular, properties of the structure of the {\em Auslander-Reiten quiver} of the derived category $\Der(A)$. Indeed, let $\Gamma_A$ be the Auslander-Reiten quiver of the derived category of an algebra $A$. We consider the {\em class quiver} $[\Gamma_A]$ of $A$ formed by the quotients $[{\mathcal C}]$ of components $\mathcal C$ of $\Gamma_A$, obtained by identifying $X, Y \in {\mathcal C}$ if $[X]=[Y]$ in $K_0(\Der(A))$. Moreover, we recall that a triangular algebra $A$ is said to be $\frac{p}{q}$-{\em Calabi-Yau} for integers $q \ge 1$ and $p \in \Z$ if $S^q=[p]$ in the derived category $\Der(A)$, where $S=\tau \circ[1]$. 
\vskip.4cm
\begin{theorem}\cite{Pe} \quad
Let $A$ be a triangular algebra with $n$ non-isomorphic simple modules and $\chi_A=\prod_{m \in M}\Phi_m^{e(m)}$ be an irreducible decomposition of its Coxeter polynomial. Let $\Gamma_A$ be the Auslander-Reiten quiver of the derived category $\Der(A)$. 
Then the following holds: 

(a) if $A$ satisfies the Calabi-Yau property then $\phi_A$ is periodic;

(b) every component of the class quiver $[\Gamma_A]$ is a tube if and only if  $\phi_A$ is periodic. 

If $\phi_A$ is periodic then $A$ is of cyclotomic type and 

(c) the period is ${\rm l.c.m.}\{\phi(m) : m \in M\}$, where $\phi(-)$ is Euler's totient function;

(d) every component of $\Gamma$ is either a tube $\Z \A_\infty / (p)$ of finite period $p$ where $p=\phi(m)$ for some $m \in M$ or of the form $\Z \Delta$ for $\Delta$ a Dynkin or extended Dynkin diagram or of one of the shapes $\A_\infty$, $\A_\infty^\infty$ or $\D_\infty$. Moreover,

(e) if $p_1,\ldots,p_s$ are the periods $> 1$ of non-homogeneous tubes in $[\Gamma]$ then $\sum\limits_{i=1}^s p_i \le n$.
\end{theorem}
\vskip.5cm

\section{Mahler measure}
Given a Laurent polynomial $P$ with integer coefficients, its Mahler's measure $M(P)$ is defined as the
geometric mean of the function $|P|$ over the real circle, that is,
$$M(P) = {\rm exp} \left(\int_0^1 {\rm ln} (|P(e^{2 \pi i t})|) dt \right)$$
\noindent
see \cite{Ma}. For a polynomial in one variable, Mahler obtained the more  elementary expression for the measure given in the Introduction.

\subsection{}
Let $P= a_0+a_1 T+ \cdots +a_n T^n$ be any polynomial with real coefficients. Besides the Mahler measure of $P$, we shall consider other measures, such as,
$$L(P)= \sum_{i=1}^n |a_i|, \,\, ||P||= {\rm max}_{|x|=1} |P(x)|$$
\indent
There are the following basic inequalities, shown in \cite{Ma} (see also \cite{Pr})
$$M(P) \le ||P|| \le L(P) \le 2^n \, M(P). $$
\noindent

\subsection{}
Let $A$ be an algebra and $\lambda_1, \ldots, \lambda_n$ be the eigenvalues  of the Coxeter transformation $\phi_A$.
We may assume that $|\lambda_1| \le |\lambda_2| \le \ldots \le |\lambda_n|=\rho_A$.
Part of the next statement was shown in \cite{Pe}.

\begin{propo}
The following inequatilities hold:
$$1 \le M(\chi_A)^{1/n} \le \frac{1}{n}e(\phi_A) \le \frac{1}{n}\, \sum\limits_{i=1}^n |\lambda_i|^2 \le \rho_A^2 \le M(\chi_A)^2$$
\end{propo}
\begin{proof}
The three inequalities not involving the Mahler measure were proved in \cite{Pe}, Proposition (3.2). The upper bound for the Mahler measure follows from the arithmetic mean-geometric mean inequality for the set of absolute values
 $\{ |\lambda_1|,\ldots,|\lambda_n| \}$. The lower bound $\rho_A \le M(\chi_A)$ is evident.
\end{proof}

\section{Accessible algebras}

\subsection{}
For an algebra $B$ we say that $A$ is {\em accessible from\/} $B$ if there is a sequence $B=B_1,B_2,\ldots,B_s=A$ of algebras such that each $B_{i+1}$ is a one-point extension (resp.\ coextension) of $B_i$ for some exceptional $B_i$-module $M_i$.

As a special case, a $k$-algebra $A$ is called {\em accessible} if $A$ is accessible from $k$. By construction, each accessible algebra $A$ is connected, moreover the indecomposable projective  $A$-modules can be arranged to form  an exceptional sequence $(P_1,\ldots,P_n)$, that is, $\Ext_A^i(P_j,P_s)=0$ for every $j > s$ and $i \ge 0$. In particular, the quiver of an accessible algebra $A$ has no loops or oriented cycles (we say that $A$ is {\em triangular}) and therefore $A$ has finite global dimension.
If not stated otherwise we shall assume that an algebra $A$ is defined over $k$ and {\em connected}, that is, its quiver is connected. Many well-known examples of algebras are accessible: hereditary algebras of tree type, more generally tree algebras, canonical algebras with three weights, poset algebras without crowns, representation-finite algebras with vanishing first Hochschild cohomology of global dimension $\leq2$. 

The construction of sequences of algebras $A=A_n,A_{n-1},\ldots,A_1=k$ where each algebra $A_{i+1}$ is a one-point extension or coextension of $A_i$ by an exceptional module $M_i$, for $i=1, \ldots,n-1$ allows to perform informative inductive procedures. For instance, every accessible algebra $A$ has vanishing Hochschild cohomology $H^i(A)=0$, for $i>0$ and $H^0(A)=k$. For any tree algebra $A$ such an inductive procedure yields that the Coxeter polynomial takes values $0$ or $1$ when evaluated at $-1$, see \cite{LePe3, LePe4}.

\subsection{} Assume that $A$ is a finite dimensional algebra of finite global dimension,
implying that the bounded derived category $\Der(A)$ of finite dimensional $A$-modules is homologically finite. Following \cite{Ke}, a
triangulated category ${\mathcal T}$ is called {\em homologically finite} if for any two objects $X$ and $Y$ from ${\mathcal T}$ the space $\Hom_{\mathcal T}(X,Y[n])$ is non-zero only for finitely many $n$. Note that a module $E$ is exceptional in $\mod{A}$ if and only if it is exceptional as an object in the triangulated category
$\Der(A)$ under the standard embedding from $\mod A$ to $\Der(A)$.

Consider an exceptional object $E$ in a triangulated category ${\mathcal T}$, then the {\em right perpendicular category} $E^{\perp}$  of $E$ consists of all objects $X$ from ${\mathcal T}$ satisfying the conditions  $\Hom_{\mathcal T}(E,X[n])=0$ for each integer $n$. Viewed as a full subcategory of ${\mathcal T}$, $E^{\perp}$ is a triangulated category, and the exact inclusion of $E^{\perp}$ in ${\mathcal T}$ admits an exact left adjoint $e_\lambda :{\mathcal T} \to E^{\perp}$, see \cite{BoKa}.

Let $M$ be an exceptional $B$-module. In case the perpendicular category $M^{\perp}$ formed in the derived category $\Der(B)$ is of the form $\Der(C)$ for an accessible algebra $C$ we say that $M$ is {\em special}. In this case for the one-point extension algebra $A=B[M]$ we get $\chi_A=(T+1) \chi_B- T \chi_C$, a useful formula in our considerations. 

An {\em accessible tower of algebras} is a family of algebras $A=A_n,A_{n-1},\ldots,A_1=k$ where each algebra $A_{i+1}$ is a one-point extension or coextension of $A_i$ by a special exceptional module $M_i$, for $i=1, \ldots,n-1$ such that $M_i^{\perp}$, formed in the derived category, $\Der(A_i)$ is equivalent to $\Der(A_{i-1})$. 

\subsection{}
Forming one-point extensions is in some sense inverse to forming perpendicular categories (with respect to an exceptional object). See \cite[Proposition 4.11]{GeLe} and \cite[section 18.3]{Le} for a detailed discussion. As an example, we recall the following remark from \cite{LePe3}.

\begin{propo} \label{cor:dynkin}
Let $A$ be a $k$-algebra which is derived equivalent to the path algebra of a Dynkin quiver $\Delta$. Let $B$ be a one-point-extension or coextension by an indecomposable $A$-module $M$. Then $B$ is derived equivalent to the path algebra of a quiver $\Gamma$ obtained from $\Delta$ by adding a new vertex and a new arrow.
\end{propo}

\subsection{}
Assume $\T$ is a triangulated category which is algebraic in the sense of \cite{Ke}. In this paper, an exceptional object $E$ is said to be of
{\em coherent type in $\T$} if the left (resp.\ right) perpendicular category ${}^\perp E$ (resp.\ ${E}^\perp$) is equivalent to
$\Der({\coh\X})$ for some weighted projective line $\X$ and,
moreover, the left adjoint $\ell$ (resp.\ right adjoint $r$) to inclusion maps $E$ to a line bundle in
$\coh\X$. The next proposition, taken from \cite{LePe4} is the key tool to determine the shape of the derived category of an extended canonical algebra.

\begin{propo}
Let $\T$ be a triangulated category having an exceptional object of coherent type
$E$. Then there exists a tilting object
$\bar{T}$ of $\T$ whose endomorphism ring is an extended canonical algebra. Further for $A=\End \,{\bar{T}}$ the categories $\T$ 
and $\Der(\mod{\, A})$ are equivalent as triangulated categories.
\end{propo}
\begin{proof}
By \cite{GeLe} the line bundle $rE$ of $\coh\X$ extends to a tilting bundle $T$ in $\coh\X$. Now $\bar{T}=T\oplus E$ is a tilting object in $\T$, whose endomorphism ring is an extended canonical algebra. 
\end{proof}

\subsection{Example:} As an illustration of the former definitions (4.1), (4.2) and as application of the Theorems stated at the Introduction, we consider the 
following tower of algebras.

Consider the  algebra $R_{2n}$  with $2n$ vertices and whose quiver is given as
$$\xymatrix@R12pt@C10pt{1\ar[r]\ar[d]\ar@{.}[dr]&2\ar[r]\ar[d]\ar@{.}[dr]&3\ar[r]\ar[d]&\cdots \ar[r] &n \ar[d] \\ 1'\ar[r]&2'\ar[r]&3'\ar[r]&\cdots \ar[r]&n'}
$$
with all commutative relations. The corresponding Coxeter polynomial 
$$\chi_{R_{2n}}= \chi_{\A_n}\otimes \chi_{\A_2} = v_{n+1} \otimes v_3$$ 
\noindent
is a product of cyclotomic polynomials. In fact $R_{2n}=\A_n \otimes \A_2$, where $\A_s$ is the hereditary algebra associated to 
the linear quiver $1 \to 2 \to \cdots \to s$, see \cite{Pe}.

For $2m+1$ odd, we consider 
$$
\def\ci{\circ}
R_{2m+1}\quad\xymatrix@R12pt@C10pt{
 &\ci\ar@{.}[dr] \ar[d]\ar[r] & \ci\ar[d]\ar[r] &\cdots\ar[r]&\ci\ar@{.}[dr]\ar[d]\ar[r]&\ci\ar[d] \\
\ci \ar[r]       &\ci \ar[r]       &\ci\ar[r] &\cdots\ar[r]&\ci\ar[r]& \ci\\}
$$
\noindent
The sequence of algebras $R_n$ forms an accessible tower of algebras, as it is easy to verify.  The following holds for the sequence of Coxeter polynomials $\chi_{R_n}$:

(a) $\chi_{R_n}=T^n + T^{n-1} - T^3 \, \chi_{R_{n-6}} +T+1$, for all $n \ge 6$;

(b) $L(\chi_{R_n})= 4 \, m + n_0$, if $n=6 \, m +n_0$ for $0 \le n_0 \le 5$;

(c) $M(\chi_{R_n})=1$.

\begin{proof}
(a): For simplicity, we write $\chi_n := \chi_{R_n}$. We want to show
$$\chi_n=(T+1)(T^{n-1}+1)- T^3\chi_{n-6} \hskip1cm (1,n)$$ 
\noindent
We shall consider also the following formulae:
$$\chi_n=(1+T)\chi_{n-1}-T(1+T)\chi_{n-3}+T^2(1+T)\chi_{n-5}-T^3\chi_{n-6} \hskip1cm   (0, n)$$
\noindent
$$\chi_{n+6}=T^{n+6}+T^{n+5}-T^{n+3}-T^{n+2}+T^6 \chi_{n-6}-T^4-T^3+T+1 \hskip1cm  (2,n)$$
\noindent
We observe that formula $(0,n)$ holds for all $n \ge 0$ as can be easily verified from the accessibility of the tower $R_n$.
We leave this as an exercise.

Obviously,  formula $(1,n)$ implies $(2,n)$. We claim the following assertions hold:

(i) formulas $(2,n)$ and $(1,n)$ imply $(1,n+6)$

(ii) formulas $(2,n)$ hold for all $n \ge 0$.

We observe that (i) and  (ii) imply statement (a). Indeed, we proceed by induction on $n$. Clearly, the formula $(1,n)$ holds for
$n=1, 2,\ldots, 6$. Assume $(1,m)$  holds for $m \le n-1$. Then $(2,n-6)$ and $(1,n-6)$ hold and by (i), $(1,n)$ holds.

(i): Observe that $(1,n)$ can be written as:
$$ T^3\chi_{n-6}= - \chi_n +T^n +T^{n-1}-T-1$$
\noindent
which sustituted in $(2,n)$ yields:
$$\chi_{n+6}=T^{n+6}+T^{n+5}-T^{n+3}-T^{n+2}+T^3(- \chi_n+T^n+T^{n-1}+T+1) -T^4-T^3+T+1 =$$
$$=T^{n+6}+T^{n+5}-T^3 \chi_n +T+1$$
\noindent
that is, $(1,n+6)$ holds.

(ii): By induction on $n$. The assertion $(2,n)$ holds for $n=1,2,\ldots,6$ by direct verification.

Assume $(2,m)$ holds for all $m \le n-1$, we show $(2,n)$. Use formula $(0,n+6)$ to obtain:
$$ \chi_{n+6}= (1+T)\chi_{n+5}-T(1+T)\chi_{n+3}+T^2\chi_{n+1}-T^3\chi_n=$$
$$=(1+T)(T^{n+5}+T^{n+4}-T^{n+2}-T^n+T^6 \chi_{n-7}-T^4-T^3+T+1)-$$
$$-(1+T)(T^{n+4}+T^{n+3}-T^{n+1}-T^n+T^7 \chi_{n-9}-T^5-T^4+T^2+T)+$$
$$+(1+T)(T^{n+2}+T^{n+1}-T^{n-1}-T^{n-2}+T^8 \chi_{n-11}-T^6-T^5+T^3+T^2)-$$
$$-T^{n+3}+T^{n+2}-T^{n}-T^{n-1}+T^6 \chi_{n-12}-T^7-T^6+T^4+T^3
$$
An elementary simplification and formula $(0,n-6)$ yield formula $(2,n)$. The proof of (a) is complete.

(b):  follows from the explicit calculation of $\chi_{R_n}$ for $0 \le n \le 5$.

(c): observe that $R_{2\, m}$ is of cyclotomic type for all $m$. Then (c) follows as a consequence of Theorem 2.
\end{proof}

\subsection{Extended canonical algebras of critical type.}

In \cite{LePe4} the extended canonical algebras $A$ with minimal weigth type $(p_1,p_2,\ldots, p_t)$ such that
${\rm Spec \,}(\chi_A(T)) \not\subset \s^1$ were classified. Recall that such an extended canonical algebra satisfies
$\dim_k H^2(A)= t-3$. In particular, if $A$ is accessible then $t=3$. This remark restricts the possible critical weights.

\begin{theor}
Let $A$ be an accessible algebra derived equivalent to an extended canonical algebra of
weight ${\bf p}= (p_1,p_2,p_3)$. Then the following holds:

(a) if ${\bf p}$ is minimal such that $A$ is not of cyclotomic type then $A$ belongs to {\rm Table 1}.

(b) if $A$ is not derived equivalent to a canonical algebra but $A$ is of cyclotomic type
then $A$ belongs to {\rm Table 2}. 

(c) if $A$ is in {\rm Table 1} then $M(\chi_A)=\rho_A^2$.
\end{theor}

\vskip.11cm

{\tiny
\begin{center}
\begin{tabular}{c|c|l|c|c}
&Weight &$\!\!$Irreducible factorization of $\chi_A$
  &$\!\!\rho (\varphi_A)\!\!$ &$\!\!\!$Dynkin index$\!\!\!$\\[8pt]
\hline &&&&\\[-3pt]
  &$(2,3,11)$
    &$\!\!1\!+\!T\!-\!T^3\!-\!T^4\!+\!T^6\!+\!T^7\!+\!T^9\!+\!T^{10}\!-\!T^{12}\!-\!T^{13}\!+\!T^{15}\!+\!T^{16}$
    &$\!\!1.1064\!\!$ &$6$\\[8pt]
  &$(2,4,9)$
    &$\!\!\Phi_2\Phi_5(T^{10}\!-\!T^9\!+\!T^5\!+\!T\!+\!1)$
    &$\!\!1.1329\!\!$ &$4$\\[8pt]
  &$(2,5,8)$
    &$\!\!1\!+\!T\!+\!T^4\!+\!T^5\!+\!T^6\!+\!2T^8\!+\!T^9\!+\!T^{10}\!+\!T^{11}\!+\!T^{14}\!+\!T^{15}$
    &$\!\!1.1574\!\!$ &$4$\\[8pt]
  &$(2,6,7)$
    &$\!\!1\!+\!T\!+\!T^4\!+\!2T^5\!+\!2T^6\!+\!T^7\!+\!T^8\!+\!2T^9\!+\!2T^{10}\!+\!T^{11}\!+\!T^{14}\!+\!T^{15}$
    &$\!\!1.1669\!\!$ &$4$\\[8pt]
  &$(3,3,8)$
    &$\!\!1\!+\!T\!+\!T^2\!+\!T^5\!+\!2T^6\!+\!3T^7\!+\!2T^8\!+\!T^9\!+\!T^{12}\!+\!T^{13}\!+\!T^{14}$
    &$\!\!1.1498\!\!$ &$3$\\[8pt]
  &$(3,4,7)$
    &$\!\!1\!+\!T\!+\!T^2\!+\!T^3\!+\!T^4\!+\!2T^5\!+\!3T^6\!+\!3T^7\!+\!3T^8\!+\!2T^9\!+\!T^{10}\!+\!T^{11}\!+\!T^{12}\!+\!T^{13}\!+\!T^{14}\!\!$
    &$\!\!1.1847\!\!$ &$3$\\[8pt]
  &$(3,5,6)$
    &$\!\!\Phi_3(T^{12}\!+\!T^9\!+\!T^8\!+\!T^7\!+\!T^6\!+\!T^5\!+\!T^4\!+\!T^3\!+\!1)$
    &$\!\!1.1966\!\!$ &$3$\\[8pt]
  &$(4,4,6)$
    &$\!\!\Phi_2\Phi_4(T^{10}-T^9+T^8+T^6+T^4+T^2-T+1)$
    &$\!\!1.2715\!\!$  &$3$\\[8pt]
  &$(4,5,5)$
    &$\!\!\Phi_5(T^{10}\!+\!T^7\!+\!T^6\!+\!T^5\!+\!T^4\!+\!T^3\!+\!1)$
    &$\!\!1.2277\!\!$ &$3$\\[8pt]
\hline &&&&\\[-3pt]
\end{tabular}
\vskip12pt
{\footnotesize{\bf Table 1.} Critical weight sequences.} 
\end{center}
}

\begin{center}
{\tiny
\begin{tabular}{rc|l|l|c}
&Weight sequence\ \ &\ \ Factorization of $\chi_A$\ \ &\ \ Period of $\phi_A$\ \ \\[3.5pt]
\hline
&&&&\\[-9pt]
&$(2,3,7)$ &$\Phi_{42}$ &$42$\\[3.5pt]
&$(2,3,8)$ &$\Phi_2\cdot \Phi_{10}\cdot \Phi_{30}$
   &$30$\\[3.5pt]
&$(2,3,9)$ &$\Phi_3\cdot \Phi_{12}\cdot \Phi_{24}$
   &$24$\\[3.5pt]
&$(2,3,10)$ &$\Phi_2\cdot \Phi_{16}\cdot \Phi_{18}$
  &$72$\\
\multicolumn{5}{c}{\phantom{\tiny .}\punteada\phantom{\tiny .}}\\
&$(2,4,5)$ &$\Phi_2\cdot \Phi_6\cdot \Phi_{30}$
  &$30$\\[3.5pt]
&$(2,4,6)$ &$\Phi^2_2\cdot \Phi_{22}$
  &$22$\\[3.5pt]
&$(2,4,7)$ &$\Phi_2\cdot \Phi_9\cdot \Phi_{18}$
  &$18$\\[3.5pt]
&$(2,4,8)$ &$\Phi^2_2\cdot \Phi_4\cdot \Phi_{12}\cdot \Phi_{14}$
  &$84$\\
\multicolumn{5}{c}{\phantom{\tiny .}\punteada\phantom{\tiny .}}\\
 &$(2,5,5)$ &$\Phi_5\cdot \Phi_{20}$
  &$20$\\[3.5pt]
&$(2,5,6)$ &$\Phi_2\cdot \Phi_8\cdot \Phi_{16}$
  &$16$\\[3.5pt]
&$(2,5,7)$ &$\Phi_{11}\cdot \Phi_{12}$
  &$132$\\[3.5pt]
&$(2,6,6)$ &$\Phi^2_2\cdot \Phi_3\cdot \Phi_6\cdot \Phi_{10}\cdot \Phi_{12}$ &$60$\\
\multicolumn{5}{c}{\phantom{\tiny .}\punteada\phantom{\tiny .}}\\
 &$(3,3,4)$ &$\Phi_3\cdot \Phi_{24}$
  &$24$\\[3.5pt]
&$(3,3,5)$ &$\Phi_2\cdot\Phi_3\cdot\Phi_6\cdot\Phi_{18}$
  &$18$\\[3.5pt]
&$(3,3,6)$ &$\Phi^2_3\cdot\Phi_{15}$
   &$15$\\[3.5pt]
&$(3,3,7)$ &$\Phi_2\cdot\Phi_3\cdot\Phi_4\cdot\Phi_{10}\cdot\Phi_{12}$
  &$60$\\[3.5pt]
&$(3,4,4)$ &$\Phi_2\cdot\Phi_4\cdot\Phi_{16}$
  &$16$\\[3.5pt]
&$(3,4,5)$ &$\Phi_{13}$ &$13$\\[3.5pt]
&$(3,4,6)$ &$\Phi_2\cdot\Phi_3\cdot\Phi_9\cdot\Phi_{10}$
  &$90$\\[3.5pt]
&$(3,5,5)$ &$\Phi_2\cdot\Phi_5\cdot\Phi_8\cdot\Phi_{10}$
   &$40$\\[3.5pt]
&$(4,4,4)$ &$\Phi^2_2\cdot\Phi^2_4\cdot\Phi_6\cdot\Phi_{12}$
  &$12$\\[3.5pt]
&$(4,4,5)$ &$\Phi_2\cdot\Phi_4\cdot\Phi_8\cdot\Phi_9$
  &$72$\\
\multicolumn{5}{c}{\phantom{\tiny .}\punteada\phantom{\tiny .}}\\
\end{tabular}
}
\vskip6pt
{\footnotesize {\bf Table 2.} Weights $p$ with $M(\phi_A)=1$}.
\end{center}

\subsection{Proof of Theorem 1.} In fact, we show the next, more precise, statement.

\begin{theor} Let $A$ be an accessible algebra not of cyclotomic type. Then there is a convex subcategory $B$ of $A$ satisfying the following properties:

(a) $B$ is minimal not of cyclotomic type, that is, if $C$ is any proper convex subcategory of $B$,  then $C$ is of cyclotomic type;

(b) $B$ is derived equivalent to a minimal wild hereditary algebra;

(c) the Mahler measure of $B$ is $M(\chi_B) \ge \mu_0$.
\end{theor}
\begin{proof}
Let $k=A_0, A_1, \ldots,A_n=A$ be a tower algebras such that $A_{i+1}$ is a one-point extension or coextension of $A_i$ by an exceptional $A_i$-module $M_i$. Since $A$ is not of cyclotomic type, there is a first index $2 \le j \le n$ such that
$A_j$ is not of cyclotomic type. We can choose that chain in such a way that $A_j$ has minimal possible $k$-dimension. The selected algebra $B=A_j$ 
clearly satisfies (a). 

We observe that algebras in Table 1 contain a convex subcategory which is a minimal wild hereditary algebra of type $[2,4,5]$ or $[2,3,7]$. Therefore to show (b) and (c)
it is enough to show that $B$ is either wild hereditary or belongs to Table 1. 

Write $B$ as a one-point extension $B=C[M]$ of an algebra $C$ of cyclotomic type by an exceptional $C$-module $M$. We distinguish two situations:

(1) $B$ is piecewise hereditary. Then also the algebra $C$ is derived equivalent to a hereditary category ${\mathcal H}$. 
Then, invoking \cite{Ha1}, ${\mathcal H}$ is either derived equivalent to the module category of a hereditary algebra $H$ or to the category of coherent sheaves $\coh \X$ over a weighted projective line $\X$. In either case $\Der({\mathcal H})$ is equivalent to the repetitive category $\bigvee_{n \in \Z}{\mathcal H}[n]$ of $\mathcal{H}$. Hence we can assume that $M \in {\mathcal H}$ and $M_{\Der({\mathcal H})}^{\perp}=\Der(M_{\mathcal{H}}^{\perp})$. We distinguish two subcases:

$\bullet$ if ${\mathcal H}$ is derived equivalent to the module category of the hereditary algebra $H$, then $M_H^{\perp}=\mod(D)$ for some hereditary algebra $D$. Since $C$ is of cyclotomic type then $H$ is hereditary of extended Dynkin type ${\tilde \Delta}$ and  $D$ is of Dynkin type. Therefore $M$ is an exceptional $C$-module and $B=C[M]$ is derived equivalent to a wild hereditary algebra $kQ$, where $Q$ is a quiver with vertices those of ${\tilde \Delta}$ and an additional source. Moreover, since $B$ is minimal not of cyclotomic type then $Q$ is of double extended Dynkin type ${\tilde {\tilde \Delta}}$ and by 
(2.5) we have $M(\chi_B) \ge \mu_0$.

$\bullet$ If ${\mathcal H}$ is derived equivalent to the coherent sheaves over $\X$, either $M$ is a vector bundle or a finite length module. In the first case, observe that, by property (a), $C$ is canonical of extended Dynkin type or a tubular algebra,
of extension type $p_0=(p_{01},p_{02},p_{03})$ for some numbers $p_{01} \le p_{02}\le p_{03}$.
Consider a filtration 
$0=M_0 \subset M_1 \subset M_2 \subset \ldots \subset M_r=M$ such that $L_i = M_i/M_{i-1}$
is a line bundle, for $i=1,\ldots, r$. 
Therefore $L_1$ is exceptional of coherent type. Either $r=1$ and, as above, $B=C[L_1]$ is an extended canonical algebra 
in {\bf Table 1}, or $r>1$ and, by property (a), $B_1=C[L_1]$ is an extended canonical algebra in {\bf Table 2}. Since the
extension type of $B_1$ is $p_1=(p_{01},p_{02},p_{03}+1)$ then, property (a) again, yields that $B_1$ is canonical of extended Dynkin type or a tubular algebra. Then $B_2=B_1[L_2]=C[M_2]$ is an extended canonical algebra. Either $r=2$ and 
$B_2$ is an extended canonical algebra in {\bf Table 1}, or $r>2$ and $B_2$ is an extended canonical algebra in {\bf Table 2}.
We may continue to reach $B=C[M]$ as an extended canonical algebra in {\bf Table 1} and $M(\chi_B) \ge \mu_0$.

Finally, in case $M$ is a finite length module of regular length $m$, $M$ lies in a tube ${\mathcal T}$ of rank $r >m$ since $M$ is exceptional. Then $M_{\mathcal H}^{\perp}= \mod(H') \times \coh(\X')$, where $H'$ is a hereditary algebra of linear type $\A_m$ and $\X'$ is a weighted projective line (of a weight dominated by that of $\X$). Therefore  $M_{\Der({\mathcal H})}^{\perp}=\Der(H') \times \Der(C')$ for a canonical algebra $C'$ whose module category is derived equivalent to $\coh(\X')$. Therefore $M$ is special exceptional and, as above, $B$ belongs to {\bf Table 1} and $M(\chi_B) \ge \mu_0$.

(2) $B$ is not piecewise hereditary. Let $D$ be  maximal convex subcategory  of $B$ such that $D$ is piecewise hereditary. Then $D$ is of cyclotomic type. Assume $D$ is derived equivalent to the hereditary category ${\mathcal H}$ and there is an exceptional $D$-module $N$ such that $D[N]$ (or equivalently, $[N]D$) is a convex subcategory of $B$ which is not derived hereditary.  We distinguish two situations:

$\bullet$ if ${\mathcal H}$ is derived equivalent to the module category of the hereditary algebra $H=k\Delta$. 
We consider the different representation types of $H$.

(i) Assume that $\Delta$ is of Dynkin type. By the Proposition above, the one-point extension $D[N]$ is derived hereditary, a contradiction to the maximality of $D$.

(ii) Assume that $\Delta$ is of extended Dynkin type. Then $\mod H$ is derived equivalent to the coherent sheaves over a weighted projective line, a situation that we shall consider later.

(iii) Assume $H$ is a wild hereditary algebra. By \cite{Xi}, the spectral radius $\rho_H \ge \mu_0$ and therefore $D$ is not of cyclotomic type, a contradiction to property (a) satisfied by $B$.

$\bullet$ If ${\mathcal H}$ is derived equivalent to the coherent sheaves over $\X$, either $N$ is a vector bundle 
or a finite length module. In the first case, $N_{\mathcal H}^{\perp}=\mod(F)$ for a hereditary algebra $F$ and $D$ is a one-point extension of $F$ by an exceptional module $L$, say $D=F[L]$. As above, $F$ is of extended Dynkin type ${\tilde \Delta}$.  
Let $Q$ be a quiver of type $\tilde{\Delta}$, therefore admitting a unique
positive {\em additive function} $\lambda$ assuming the value $1$
at some  extension vertex $p$. Additivity of $\lambda$ means that for any vertex $u$ one has 
$2\lambda(u)=\sum_{v }a_{uv}\lambda(v)$, where $v$ runs through all vertices and $a_{uv}$ denotes the number of edges between $u$ and $v$.
Then $D$ is derived equivalent to the path algebra $\tilde{\tilde{\Delta}}$ which is the graph arising from $\tilde{\Delta}$
by adjoining a new edge at the extension vertex $p$, see \cite{LePe3}. Therefore $D$ is hereditary of wild type then not cyclotomic, again a contradiction.

In case $N$ is a finite length module of regular length $m$, $N$ lies in a tube ${\mathcal T}$ of rank 
$r >m$ since $N$ is exceptional. Then $N_{\mathcal H}^{\perp}= \mod(H') \times \coh(\X')$, where $H'$ is a hereditary algebra of linear type $\A_m$ and $\X'$ is a weighted projective line of a weight dominated by that of $\X$. This contradicts
that $F$ is of extended Dynkin type and hence the weight of $\X$ is minimal. In conclusion, case (2) never arises.
\end{proof}

\section{Interlaced towers of algebras}

\subsection{}
Following \cite{LePe2}, we shall say that a
polynomial $f \in \Z [T]$ is {\em represented\/} by $q \in \Z
[T]$ if $f(T^2)=q^*(T):=T^{\deg(q)}q(T+T^{-1})$. 
Two important families of representable polynomials are the following:

$\bullet$ A fundamental fact for a hereditary algebra $A=kQ$, when $Q$ is a
{\em bipartite quiver\/} without oriented cycles is that
$\Spec\,\varphi_A\subset \s^1\cup \R^+$. This was shown by A'Campo
\cite{AC} as a consequence of the identity $\chi_A(T^2)=T^n\kappa_{|Q|}
(T+T^{-1})$, where $n$ is the number of vertices of $Q$ and
$\kappa_{|Q|}$ is the characteristic polynomial of the adjacency
matrix of the underlying graph $|Q|$ of $Q$. 

$\bullet$ Recall  that the (normalized) {\em Chebycheff polynomials} (of
the second kind) $(u_n)_n$ may be inductively constructed by the rules:
$$ u_0=1,\; \ u_1=T, \; \textrm{ and } \;
u_{n+1}=Tu_n-u_{n-1}, \textrm{ for }n\ge 1.$$
A simple induction shows that the characteristic polynomial of the Coxeter transformation of a linear graph
$\A_n = [n]$ is the  polynomial $u_n(T)$. The  polynomial $u_n$ represents the polynomial $v_{n+1}=T^n+T^{n-1}+\cdots+T+1$.

Other examples of representable polynomials are the cyclotomic polynomials and
the Coxeter polynomials of canonical and extended canonical algebras, as shown in \cite{LePe2}. We shall come back to some of these results.

The following basic remark was proved in \cite{LePe2}. We denote $\Root(p)$ the set of (complex) roots of a polynomial $p$.

\begin{propo}
Let $p$ and $q \in \Z[T]$ be such that $p$ is represented by $q$. Assume that $p(T)$ is a self-reciprocal polynomial. Then

\begin{itemize}
\item[{\rm (a)}] $\mu \in \Root(p)$, then
$\mu^{-1},\bar{\mu},\bar{\mu}^{-1}\in \Root(p)$;

\item[{\rm (b)}] $\Root(p)\subset \s^1\cup \R^+$ if and only if $\Root(q)\subset \R$;

\item[{\rm (c)}] $\Root(p)\subset \s^1$ (resp. $\s^1\setminus
\{1\}$) if and only if $\Root(q)\subset [-2,2]$ (resp. $(-2,2)$);

\item[{\rm (d)}] Assume $q(T)=\chi_\Delta (T)$ is the characteristic polynomial of a the adjacency matrix of a graph 
$\Delta$ then

\begin{itemize}
\item[{\rm (i)}] $\Delta$ is a bipartite graph.

\item[{\rm (ii)}] $\Root(q)\subset \R$ and $\Root(p)\subset \s^1\cup \R^+$.

\item[{\rm (iii)}] $\Root(p)\subset \s^1\setminus \{1\}$ if and only if $\Delta$ is a union of Dynkin graphs.
\end{itemize}
\end{itemize}
\end{propo}

We shall study representability of Coxeter polynomials associated to towers of algebras and obtain consequences.

\subsection{}
Let $B=A_m, A_{m+1}, \ldots,A_n=A$ be a family of triangular algebras with $m \le n-2$ satisfying, for each $s=m,\ldots,n$:

(i) each algebra $A_s$ has $s$ simple pairwise non-isomorphic modules;

(ii) the Coxeter polynomial $\chi_{A_s}$ is represented by a monic polynomial $p_s$ of degree $s$;  

(iii) $p_s=T p_{s-1} - p_{s-2}$, if $s \ge m+2$;

(iv) $p_m$ has $m$ pairwise different real roots. 

If conditions (i) to (iv) are satisfied, we say that $B=A_m, A_{m+1}, \ldots,A_n=A$ is an {{\em interlaced tower of algebras}.
The next statement is the motivation for the introduction in \cite{LePe2} of the representability concept. We shall recall here the simple proof. 

\begin{propo}
Let $B=A_m, A_{m+1}, \ldots,A_n=A$ be a family of triangular algebras with $m \le n-2$ satisfying, for each $s=m,\ldots,n$, conditions {\rm (i)} and {\rm (ii)}.
Then the following hold:

${\rm (a) \,}$ Condition {\rm (iii)} is equivalent to:
$\chi_{A_{s+1}} =(T+1) \chi_{A_s} - T \chi_{A_{s-1}}.$

${\rm (b) \,}$ if moreover, conditions {\rm (iii)} and {\rm (iv)} hold, the next conditions are also satisfied:

\begin{itemize}
\item[{\rm (v)}] if two successive functions among
$p_{s-1}$, $p_s$ and $p_{s+1}$ have a common root $\lambda$ then $p_m(\lambda)=0$;

\item[{\rm (vi)}] if one of the functions vanishes in $\lambda
\in \R$, then the immediate successor and predecessor of that
function take real values of different sign at $\lambda$ (that is,
if $p_s(\lambda)=0$, then both
$p_{s-1}(\lambda)$ and $p_{s+1}(\lambda)$ are real and
$p_{s-1}(\lambda)p_{s+1}(\lambda)<0$).
\end{itemize}
\end{propo}
\begin{proof}
Assume that $\chi_s(T^2)=T^s p_s(T+T^{-1})$ for the polynomial $p_s$ representing  $\chi_s=\chi_{A_s},\, s=m,\ldots,n$.
Assume condition (iii) holds. Then
$$\chi_{s+1}(T^2)=T^{s+1}p_{s+1}(T+T^{-1})=T^{s+1}\{(T+T^{-1})p_s(T+T^{-1})-p_{s-1}(T+T^{-1})\}=$$
$$=(T^2+1)T^s p_s(T+T^{-1})-T^2T^{s-1}p_{s-1}(T+T^{-1})=(T^2+1)\chi_s(T^2)-T^2 \chi_{s-1}(T^2)$$
\noindent
that is, $\chi_{s+1}(T)=(T+1)\chi_s(T)-T\chi_{s-1}(T)$. To prove the converse, follow the steps in the reverse order.

(b) follows easily.
\end{proof}

\subsection{Examples:}

(1) Let $L$ be a connected convex subcategory of the hereditary quiver algebra $H=k\Delta$ of tree type. Then there exists a (non-unique) tower of algebras
$H_m=L,H_{m+1},\ldots,H_n=H$ such that $H_s$ is a one-point extension or coextension]1W of $H_{s-1}$ by an indecomposable module $M_{s-1}$. Clearly,  the modules $M_{s-1}$ are projective or injective, therefore excepcional. Hence $H_s=k\Delta_s$ where $\Delta_{s-1}$ is obtained from $\Delta_s$ by deleting a vertex $a_s$ and the perpendicular category $M_s^\perp$ formed  $\Der(H_s)$ is derived equivalent to $\mod(H_{s-1})$, see \cite{LePe3, LePe4} for details. Moreover, the result quoted in (5.1) yields that the Coxeter polynomial $\chi_{H_s}$ is represented by $\kappa_{|\Delta_s|}$ the characteristic polynomial of the adjacency matrix of the underlying graph $|\Delta_s|$ of $\Delta_s$. 
In case $m \le n-2$, these remarks show:

$\bullet$ the sequence $H_m=L,H_{m+1},\ldots,H_n=H$ is an accessible tower of algebras;

$\bullet$ the sequence $H_m=L,H_{m+1},\ldots,H_n=H$ is an interlaced tower of algebras.

(2) Consider the accessible tower of algebras $(R_n)$ as given in example (4.6). We observe that $\chi_n=\chi_{R_n}$ can be represented in the following way:

{\rm (i)} For $n \le 6$ the Coxeter polynomials of the sequence are:
$$\chi_0=1, \chi_1=T+1, \chi_2=T^2+T+1, \chi_3=T^3+T^2+T+1,$$
$$\chi_4=T^4+T^3+T+1,\chi_5=T^5+T^4+T+1, \chi_6=T^6+T^5-T^3+T+1$$
\noindent
represented by the polynomials
$$r_0=1, r_1=T, r_2=T^2-1,r_3=T^3-2T, r_4=T^4-3T^2,$$
$$ r_5=T^5-4T^3+2T,r_6=T^6-5T^4+5T^2-1$$

{\rm (ii)} Observe that $(T^{n-1} +1)=v_n-Tv_{n-2}$ then $T^n+T^{n-1}+T+1=(T+1)(T^{n-1}+1)$ is represented by
$w_n=T(u_{n-1}-u_{n-3})$.

{\rm (iii)} $n \ge 6$: we define  $r_n=w_n-T^3 r_{n-6}$. We verify by induction on $n$ that $r_n$ represents $\chi_{R_n}$:
 $$\chi_{R_n}(T^2)= (T^2+1)(T^{2n-2}+1)- T^6 \, \chi_{R_{n-6}}(T^2) = $$
$$=T^n \, w_n(T+T^{-1})-T^6 T^{n-6}\, r_{n-6}(T+T^{-1})=T^n \, r_n(T+T^{-1})$$

{\bf Claim:} The tower $(R_n)$ is interlaced. In fact, it is enough to show that $w_n=Tw_{n-1}-w_{n-2}$. Indeed,
$$Tw_{n-1}-w_{n-2}=T\{T(u_{n-2}-u_{n-4})\} - T(u_{n-3}-u_{n-5})=$$
$$=T(Tu_{n-2}-u_{n-3})-T(Tu_{n-4}-u_{n-5})=T(u_{n-1}-u_{n-3})=w_n$$

(3) Let $H$ be the star of type $[p_1,p_2,p_3]$ with Coxeter
polynomial $\chi_{[p_1,p_2,p_3]}$ represented by $\kappa_{[p_1,p_2,p_3]}$. 
Observe that the corresponding canonical algebra $C_{(p_1,p_2,p_3)}$ is essentially unique. 
The Coxeter polynomial $\chi_{(p_1,p_2,p_3)}$ of the canonical algebra satisfies
$$\chi_{(p_1,p_2,p_3)}=\chi_{K_2}\prod\limits^t_{i=1}\kappa_{[p_i-1]}$$
where
$\kappa_{[p]}$ denotes the characteristic polynomial of the underying graph of the linear
path $[p]$ and $K_2$ is the Kronecker diagram $\xymatrix{\bulito
\ar@{-}@/^/[r]&\bulito \ar@{-}@/^/[l]}$.

Since $\chi_{K_2}$ is represented by $(T-1)^2$ and $\kappa_{[p_i-1]}=v_{p_i}$ is represented by $u_{p_i}$, then
$\chi_{(p_1,p_2,p_3)}$ is represented by $u_{(p_1,p_2,p_3)}=(T-1)^2\prod_{i=1,2,3} u_{p_i}$.

For each number $3 \le n \in \N$ fix a triple $p(n)=(p_1(n),p_2(n),p_3(n))$ such that the following are satisfied:

$\bullet$ set $p(3)=(1,1,1)$;

$\bullet$ $p_1(n)+p_2(n)+p_3 (n)=n$;

$\bullet$ $p(n) \le p(n+1)$ coordinate wise.

{\bf Claim:} The tower of canonical algebras $C_{p(n)}$, for $3 \le n \in \N$, is interlaced. Indeed, it is enough to show that $Tu_{p(n)}-u_{p(n-1)}=u_{p(n+1)}$. We have to distinguish two cases:

$\bullet$ $p_1(n+1)=p_1(n)+1=p_1(n-1)+2$; 

$\bullet$ $p_1(n+1)=p_1(n)+1$ and $p_2(n)=p_2(n-1)+1$.

We only analize the second case, namely,

$$Tu_{p(n)}-u_{p(n-1)}=T(T-1)^2 u_{p_1} u_{p_2+1}u_{p_3}-(T-1)^2\prod_{i=1,2,3} u_{p_i}=$$
$$=(T-1)^2u_{p1+1}u_{p_2+1}u_{p_3}=u_{p(n+1)}$$

(4) Let $E=E(p_1,p_2,p_3)=C_{(p_1,p_2,p_3)}[P]$ be the extended canonical algebra of $C_{(p_1,p_2,p_3)}$ by the indecomposable projective module $P$
associated to the unique source.
Following \cite{LePe2}, the Coxeter polynomial of $E$ is $\chi_{\langle p_1,p_2,p_3 \rangle}=(T+1)\chi_{(p_1,p_2,p_3)}-T\chi_{[p_1,p_2,p_3]}$ and it is represented by 
$$q_{(p_1,p_2,p_3)}=T\chi_{K_2}\prod_{i=1,2,3}u_{p_i}-\kappa_{[p_1,p_2,p_3]}$$
\noindent
As in (3),  for each number $3 \le n \in \N$ fix a triple $p(n)=(p_1(n),p_2(n),p_3(n))$ satisfying the above conditions. As above, it is simple to show:

{\bf Claim:} The tower of extended canonical algebras $E_{p(n)}$, for $3 \le n \in \N$, is interlaced.

(5)  Consider the canonical algebra $C$ of type $(2,3,6)$ and $A=C[P]$ where $P$ is the indecomposable projective $C$-module associated
to the unique source of $C$. Then $A$ is the extended canonical algebra associated to the quiver
$$
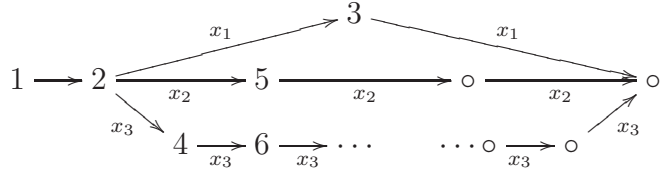
\def\c{\circ}
\xymatrix@C18pt@R10pt{
                                      & & &    &3 \ar[rrrd]^{x_1}& &&\\
1 \ar[r]&2\ar[rr]_{x_2}\ar[rrru]^{x_1}\ar[rd]_{x_3}&&5 \ar[rr]_{x_2}&  &\c\ar[rr]_{x_2} &&\c    \\
                                     &   &4\ar[r]_{x_3}          &6 \ar[r]_{x_3}          &\cdots&\cdots\c\ar[r]_{x_3}&\c\ar[ru]_{x_3}&    \\
}
$$
satisfying the equation:
$$x_3^n=x_1^2- x_2^3$$
\noindent

Consider the tower of accessible algebras $k=A_1, A_2,\dots, A_{10}$ where $A_i$ is the subcategory of $A$
generated by $1,\ldots, i$, for $i=1,\ldots,6$ and, in general, $A_{s+1}$ is a one-point coextension of $A_s$, for $s=1,\ldots,10$.

Observe that $A_6$ is hereditary of wild representation type therefore $M(\chi_{A_6})>1$. Moreover, since $\phi_{A_{10}}=\Phi_{42}$ then $M(\chi_{A_{10}})=1$.
Theorem 2 implies that the tower $(A_s)$ is not an accessible tower of algebras, neither an interlaced tower of algebras.

\subsection{} The following criterion was observed in particular cases in \cite{LePe2,LePe3,LePe5}.

\begin{propo}
Let $B=A_m, A_{m+1}, \ldots,A_n=A$ be an interlaced tower of algebras. Then the following happens:

(a) there exists an integer $p$ such that $\chi_{A_s}(-1)=0$ or $\chi_{A_s}(-1)=p^2$, for $m \le s \le n$, depending if $s$ is odd or even;

(b) $0 \le p \le 1$ if $A_m$ is derived equivalent to a tree algebra; $0 \le p \le 2$ if $A_m$ is derived equivalent to a canonical algebra or to an extended canonical algebra.
\end{propo}
\begin{proof}
It was observed in \cite{LePe3} that  if $s$ is odd $\chi_{A_s}(-1)=0$, and if $s$ is even $\chi_{A_s}(-1)=p_s^2$ for some integer $p_s$.
For any $m+1 \le s \le n-1$ we have $\chi_{A_{s+1}}=(T+1)\chi_{A_s}-T \chi_{A_{s-1}}$ which implies $\chi_{A_{s+1}}(-1)=\chi_{A_{s-1}}(-1)$. 
It follows that all $p_s$ are equal for $s$ even.

(b): this is shown in \cite{LePe2}.
\end{proof}

\subsection{} The {\em interlacing properties} of eigenvalues of Coxeter polynomials of an interlaced tower of algebras are contained in the following result 
which generalizes those of \cite{LePe2}. 

Consider an interlaced tower of algebras $B=A_m, \cdots, A_n=A$ with $m \le n-2$ and $p_s$ the polynomial representing $\chi_{A_s}$, for $s=m,\ldots,n$.
At the level of the family $p_m,\ldots,p_n$ satisfying conditions $(i)$ to $(vi)$ above, a version of {\em Sturm's
Theorem}~\cite{Ob} assures that given any interval $[\alpha,\beta]\subset \R$ and roots $\lambda_1\le \cdots \le \lambda_r$
of $p_{s+1}$ in $[\alpha,\beta]$, then $p_s$ has roots $\lambda'_1\le \cdots \le\lambda'_{r-1}$ in $[\alpha,\beta]$ satisfying $\lambda_1\le
\lambda'_1\le \lambda_2\le \lambda'_2\le \cdots \le \lambda_{r-2}\le \lambda'_{r-1}\le \lambda_r$ (interlacing property). 

\begin{theor}
Consider an interlaced tower of algebras $B=A_m, \cdots, A_n=A$ with $m \le n-2$. Then the following hold:

\begin{itemize}
\item[{\rm (a)}] If the eigenvalues of $\phi_{A_n}$ are contained in $\s^1 \cup \R^+$ (resp. in $\s^1$), then the same holds for 
$\phi_{A_s}$, for $s=m,\cdots,n$;

\item[{\rm (b)}] suppose that the eigenvalues of $\phi_{A_n}$ are contained in $\s^1 \cup \R^+$, then the following holds:

{\rm (i) \,} if $\mu_1\le \cdots \le \mu_s$ (resp. $\mu'_1\le \cdots \le \mu'_t$) are the eigenvalues of $\chi_{A_n}$ (resp. $\chi_{A_{n-1}}$) not in $\s^1$, then
$s-1 \le t \le s$ satisfying 
$$\mu_1\le \mu'_2 \le \mu_2\le \mu'_3 \le \cdots \le \mu_{s-1}\le \mu'_s \le \mu_s$$
\noindent
and possibly $1 < \mu'_1 \le \mu_1$; 

{\rm (ii) \,} either both $A_{n-1}$ and $A_n$ are of cyclotomic type or $M(\chi_{A_{n-1}}) <M(\chi_{A_n})$.

\item[{\rm (c)}] $M(\chi_{A_{m+2}}) \le M(\chi_{A_{m+1}})L(\chi_{A_m})$$
\noindent$
\end{itemize}
\end{theor}
\begin{proof}
(a): If $\chi_{A_n}$ has all its roots in $\s^1 \cup \R^+$ (resp. in $\s^1$), then 
${\rm Spec}\,(p_n) \subset \R^+ $ (resp. in $[-2,2]$). By Sturm's
Theorem, ${\rm Spec\, }(p_{n-1}) \subset \R^+ $ (resp. in $[-2,2]$) and hence
${\rm Spec\,}\chi_{A_{n-1}}\subset \s^1 \cup \R^+$ (resp. in $\s^1$).

(b) Consider the roots of $p_n$ satisfying $2 < \lambda_1<\lambda_2<\cdots <\lambda_{s-1}<\lambda_s$ and those roots
$2 < \lambda'_1<\lambda'_2<\cdots <\lambda'_t$ of $p_{n-1}$. Then either
$\lambda_1<\lambda'_1<\lambda_2<\cdots <\lambda_{s-1}<\lambda'_t<\lambda_s$
or 
$\lambda'_1<\lambda_1<\lambda'_2<\cdots <\lambda_{s-1}<\lambda'_t<\lambda_s$,
accordingly $t=s-1$ or $t=s$. 
 
Since $T^n \, p_n(T+T^{-1})= \chi_{A_n}(T^2)$ implies that $\lambda_i=\nu_i+\nu_i^{-1}$ and $\mu_i=\nu_i^2$
for any root $2 <\lambda_i$ of $p_n$, we get 
$$\nu_i={\lambda_i^2 \over 2}-1 +{\lambda_i \over 2} \sqrt{\lambda_i^2-4}$$
\noindent
Part (i) follows immediately. 

For (ii), observe that $M(\chi_{A_n})=\prod_{i=1}^s \mu_i$ and use the interlacing result just proved
to get $M(\chi_{A_{n-1}}) \le M(\chi_{A_n})$. If the equality holds then either $s=0$, implying that both $A_n$ and $A_{n-1}$ are
of cyclotomic type, or $t=s$ and $\mu'_i=\mu_i$, for all $1 \le i \le s$,
which contradicts that $p_{n-1}$ and $p_n$ have no common roots. 

(c): We have $\chi_{A_{m+2}} =(T+1) \chi_{A_{m+1}} - T \chi_{A_m}.$ Then for $x \in \s^1$ we get
$$|\chi_{A_{m+2}}(x)| \le |(x+1) \chi_{A_{m+1}} (x)|+|x \chi_{A_m}(x)| \le 2 ||\chi_{A_{m+1}}||+||\chi_{A_m}||.$$
\noindent
Therefore, 
$$M(\chi_{A_{m+2}})= {\rm exp} \left( \int_0^1 {\rm ln} (|\chi_{A_{m+2}}(e^{2 \pi i t})|) dt \right) \le $$
$$\le  {\rm exp} \left( \int_0^1 {\rm ln} (|(e^{2 \pi i t}+1) \chi_{A_{m+1}} (e^{2 \pi i t})|) dt + \int_0^1 {\rm ln} (L(\chi_{A_m})) dt \right)=$$
$$=M((T+1)\chi_{A_{m+1}})M(T\chi_{A_m})=M(\chi_{A_{m+1}})L(\chi_{A_m})$$
\noindent
where the first inequality is due to ${\rm ln}(x+y) \le {\rm ln}(x)+{\rm ln}(y)$ for $x,y \ge 0$, while the last equality is due to the multiplicative property of 
the Mahler measure.
\end{proof}

\subsection{Proof of Theorem 2:} Let $A_m,\ldots,A_n$ be an interlaced tower of algebras with $m<n$. 
Suppose that ${\rm Spec \,}\phi_{A_n} \subset \s^1\cup \R^+$. 
Then (a) of the above Theorem implies that all eigenvalues of $\phi_{A_s}$ lie in $\s^1\cup \R^+$ and therefore 
$M(\chi_{A_s})=\prod_{\mu \in {\rm Spec \,}\phi_{A_s}\setminus \s^1} \mu$, for $s \le n$.

We distinguish two cases.
Assume that $A_n$ is of cyclotomic type. Then, by induction, using (a) of the above Theorem we get that also ${\rm Spec \,}\phi_{A_m} \subset \s^1$. 
That is, $A_m$ is of cyclotomic type.

If $A_n$ is not of cyclotomic type, then (b) of Theorem yields $M(\chi_{A_{n-1}})<M(\chi_{A_n})$. By induction hypothesis, we may assume that 
$M(\chi_{A_m}) \le M(\chi_{A_{n-1}})$ which implies the desired lower inequality.
\qed

{\footnotesize
\noindent
Centro de Investigaci\'on en Matem\'aticas, A.C.\\
Guanajuato 36240 M\'exico\\
\hbox{\hskip.18cm} {\footnotesize\em e-mail:\/} jap@cimat.mx \\
and\\
Instituto de Matem\'aticas, UNAM. Cd. Universitaria, \\
M\'exico 04510 D.F.\\
\hbox{\hskip.18cm} {\footnotesize\em e-mail:\/} jap@matem.unam.mx
}

\begin{thebibliography}{10000}
\bibitem{AC}
N. A'Campo, · 
Sur les valeurs propres des transformationes de Coxeter.
{\em Inventiones mathematicae},
33(1) (1976), 61--67.

\bibitem{ASS}
  I. Assem, D. Simson and A. Skowro\'nski,
  \emph{Elements of the Representation Theory of Associative Algebras 1:
  Techniques of Representation Theory}.
  {London Mathematical Society Student Texts 65},
  Cambridge University Press, 2006.

\bibitem{BaLe} M. Barot and H. Lenzing,
One-point extensions and derived equivalence, \emph{J. Algebra} \textbf{264} (2003), 1--5.

\bibitem{BoKa}
A.I. Bondal and M.M. Kapranov. Representable functors, Serre
functors, and mutations. \emph{Math. USSR, Izv.} \textbf{35}
(1990), 519--541.

\bibitem{GeLe}
W. Geigle and H. Lenzing, A class of weighted projective curves arising in representation theory of finite dimensional algebras. In \emph{Singularities, representations of algebras, and vector bundles}. \emph{Springer Lect. Notes Math.} \textbf{1273}  1987), 265--297.

\bibitem{Ha0}
D. Happel,
Hochschild cohomology of finite dimensional algebras.
In \emph{S\'{e}minaire d'Alg\`{e}bre Paul Dubreil et Marie-Paul Malliavin, 39\`{e}me Ann\'{e}e}.
Lecture Notes in Math. \textbf{1404}, Springer-Verlag 1989, 108--126.

\bibitem{Ha2}
D. Happel,
The trace of the Coxeter matrix and the Hochschild cohomology.
{\em Linear Algebra and its Applications} 258 (1997), 169-177  

\bibitem{Ha1}
D. Happel,
A characterization of hereditary categories with tilting object,
{\em  Invent. math.} 144 (2001), 381-398

\bibitem{Ha3}
D. Happel,
The Coxeter polynomial for a one-point extension algebra.
{\em Journal of Algebra} 321 (2009), 2028-2041 

\bibitem{Ke}
B. Keller: Deriving DG categories. Ann. Sci. \'{E}cole Norm. Sup. (4) \textbf{27}(1),63–102 (1994).

\bibitem{KeOs}
D.G. Kendall and R. Osborn. 
Two simple lower bounds for Euler's function.
{\em Texas J. Sci.} 17, 1965. 

\bibitem{Kr} L. Kronecker.
{\em Zwei S\"atze \"uber Gleichungen mit ganzzahligen
Coefficienten.} J. Crelle 1857, Ouvres I, p. 105--108.

\bibitem{Le}
H. Lenzing, Coxeter transformations associated with finite-dimensional
algebras. Computational methods for representations of groups and algebras,
287--308, Progr. Math., \textbf{173}, Birkh\"{a}user, Basel,  1999.

\bibitem{LePe1}
H. Lenzing and J.A. de la Pe\~{n}a,
Wild canonical algebras.
\emph{Math. Z.} \textbf{224} (1997) 403--425.

\bibitem{LePe1.a}
H. Lenzing and J.A. de la Pe\~{n}a,
Supercanonical algebras. {\em J. of Algebra} 282 (2004), 298-348. 

\bibitem{LePe2}
H. Lenzing and J.A. de la Pe\~{n}a,
A Chebysheff recursion formula for Coxeter polynomials.
\emph{Linear Algebra and its Applications}. 430 (2009) Pages 947-956 

\bibitem{LePe3}
H. Lenzing and J.A. de la Pe\~{n}a, {Spectral analysis of finite dimensional algebras and singularities}, \emph{Trends in Representation Theory of Algebras and Related Topics}, ed. A. Skowro\'nski, EMS Publishing House, Z\"{u}rich (2008) 541--588.

\bibitem{LePe4}
H. Lenzing and J.A. de la Pe\~{n}a
Extended canonical algebras and Fuchsian singularities. 
{\em Math. Z. 268} (2011), no. 1-2, 143–167.

\bibitem{LePe5}
H. Lenzing and J.A. de la Pe\~{n}a,
{Accessible algebras}. To appear J. Algebra (2013).

\bibitem{Ma}
K. Mahler. 
Lectures on trascendental Number Theory. Lecture Notes in Mathematics {\bf 546},
Springer Verlag. Berlin (1976).

\bibitem{Mo}
R. G. Mosier
{Root neighborhoods of a polynomial}
{\em Mathematics of Computation} vol. 47, 175 (1986), 5-73

\bibitem{Mo1}
M.J. Mossinghoff, 
Polynomials with small Mahler measure. {\em Mathematics of Computation} {\bf 67} (224)
(1998), 1697–1706.

\bibitem{Mo2}
M.J. Mossinghoff,
Michael Mossinghoff's Home page.

\bibitem{Obreshkoff:1963}
N. Obreschkoff,
\emph{Verteilung und Berechnung der Nullstellen Reeller Polynome.}
VEB Deutscher Verlag der Wissenchaften, Berlin,  1963.

\bibitem{Pe}
J.A. de la Pe\~na,
Algebras whose Coxeter polynomials are product of cyclotomic polynomials. To appear {\em Algebras and Representation Theory} (2013).  

\bibitem{Pr}
V. Prasolov,
Polynomials. {\em Algorithms and Computation in Mathematics} {\bf 11} Springer Verlag. Berlin (2001).

\bibitem{Ri}
C.M. Ringel, Tame algebras and integral quadratic forms.
Lecture Notes in Mathematics, \textbf{1099}. Springer-Verlag, Berlin,  1984.

\bibitem{Sm}
C. Smyth. The Mahler measure of algebraic numbers: a survey. 
In McKee, James; Smyth, Chris. Number Theory and Polynomials. London Mathematical Society Lecture Note Series. 352. 
Cambridge University Press (2008), 322–349.

\bibitem{Sk}
A. Skowro\'nski,
Simply connected algebras in Hochschild cohomologies.
In \emph{Representations of Algebras.}
Canad. Math. Soc. Conf. Proc. \textbf{14},
Amer. Math. Soc., Providence, RI, 1993, 431--447.

\bibitem{Xi}
Ch. Xi, 
On wild hereditary algebras with small growth numbers.
{\em Communications Algebra},
18:10 (1990), 3413- 3420.
\end{thebibliography}
\end{document}